\documentclass[11pt]{article}
\usepackage{amssymb}
\usepackage{makeidx}
\usepackage{amsfonts}
\usepackage{amsmath}

\newtheorem{theorem}{Theorem}
\newtheorem{acknowledgement}[theorem]{Acknowledgement}

\newtheorem{remark}[theorem]{Remark}

\input{tcilatex}

\begin{document}

\title{Schr\"{o}dinger Equation in Phase Space and Deformation Quantization}
\author{Maurice A de Gosson \\
Universit\"{a}t Potsdam, Inst. f. Mathematik\\
Am Neuen Palais 10, D-14415 Potsdam\\
maurice.degosson@gmail.com}
\maketitle

\begin{abstract}
We justify the relevance of Torres-Vega Schr\"{o}dinger equation in \ phase
space using Stone-von Neumann's theorem, and relate it to deformation
quantization.
\end{abstract}

Received ...., revised ....



\section{Introduction}

No more than a few decades ago physicists where still very reluctant to
accept the idea of a quantum mechanics in phase space (the argument usually
invoked was that Heisenberg's uncertainty principle forbids to view points
in phase space having any physical meaning). Things have changed, and
phase-space techniques are now widely used. Roughly speaking, one can say
that the tenants of phase-space quantum mechanics belong to two groups:
those working in the beautiful and deep theory of deformation quantization
based on the work by Bayen \textit{et al}. \cite{BFFLS}, and those studying
various notions of Schr\"{o}dinger's equation in phase space; one of the
most cited approaches is that of Torres-Vega and Frederick \cite{TV1,TV2}
who, proposed a whole family of Schr\"{o}dinger equations in phase space,
whose prototype is 
\begin{equation}
i\hbar \frac{\partial \Psi }{\partial t}=H\left( \tfrac{x}{2}+i\hbar \tfrac{%
\partial }{\partial p},\tfrac{p}{2}-i\hbar \tfrac{\partial }{\partial x}%
\right) \Psi  \label{erwin2}
\end{equation}%
The aim of this paper is twofold:

\begin{itemize}
\item We will show that (\ref{erwin2}) is not only equivalent to the usual
Schr\"{o}dinger equation in \textquotedblleft configuration
space\textquotedblright\ provided that one restricts the set of solutions to
a closed subspace of $L^{2}(\mathbb{R}_{x,p}^{2})$, but that it actually
corresponds to the choice of an irreducible unitary representation of the
Heisenberg group;

\item We will examine the relationship between equation (\ref{erwin2}) and
deformation quantization; we will come to the conclusion that deformation
quantization is justified by the Stone-von Neumann theorem, and that
Torres-Vega and Frederick's theory of Schr\"{o}dinger equation in
phase-space is in fact a \textit{Doppelg\"{a}nger} of deformation
quantization.
\end{itemize}

The kernel of our argument is the following observation. Recall that Schr%
\"{o}dinger, in deriving his equation, started with a Hamiltonian function $%
H=T+V$ $(T$ the kinetic energy, $V$ a potential); elaborating on Hamilton's
optical--mechanical analogy (\cite{ICP,Jammer}) he integrated the Poincar%
\'{e}--Cartan (or: action) form%
\begin{equation}
\alpha _{H}=pdx-Hdt  \label{alfach}
\end{equation}%
in order to obtain a solution of Hamilton-Jacobi's equation for $H$. This
allowed him, by an inductive argument, to derive what we call the
time-independent Schr\"{o}dinger equation which is satisfied by a stationary
matter-wave $\psi _{0}$; in a follow-up to this paper he introduced the wave
function $\psi (x,t)=e^{-iEt/\hbar }\psi _{0}$ which is, when $\psi _{0}$ is
a solution of the time-dependent equation%
\begin{equation}
i\hbar \frac{\partial \psi }{\partial t}=H\left( x,-i\hbar \frac{\partial }{%
\partial x}\right) \psi \text{.}  \label{erwin1}
\end{equation}%
Compared to the Hamilton equations 
\begin{equation*}
\frac{dx}{dt}=\frac{\partial H}{\partial p}\text{ \ , \ }\frac{dp}{dt}=-%
\frac{\partial H}{\partial x}
\end{equation*}%
from classical mechanics, Schr\"{o}dinger's equation introduces a deep
asymmetry: the variable $p$ has disappeared altogether and has been replaced
by the operator $-i\hbar \partial /\partial x$. This asymmetry comes from
Schr\"{o}dinger's honest and totally justifiable use of the action form (\ref%
{alfach}), where the variables $p$ and $x$ play asymmetric roles. Let us now
pause and ask ourselves where the interest of the action form (\ref{alfach})
comes from. Well, it mainly comes from the fact that it is a relative
integral invariant, that is, its exterior derivative $d\alpha _{H}$ is an
absolute integral invariant. It is precisely this property that allows one
to integrate Hamilton--Jacobi's equation in terms of $\alpha _{H}$. Now,%
\begin{equation*}
d\alpha _{H}=dp\wedge dx-Hdt
\end{equation*}%
has $\alpha _{H}$ as a primitive --among infinitely many other! For
instance, every differential form%
\begin{equation*}
\alpha _{H}^{\lambda }=\lambda pdx-(1-\lambda )xdp-Hdt
\end{equation*}%
obviously satisfies 
\begin{equation*}
d\alpha _{H}^{\lambda }=dp\wedge dx-Hdt
\end{equation*}%
and is hence also a relative integral invariant. Making the particular
choice $\lambda =\frac{1}{2}$ we will denote by $\beta _{H}$ the
corresponding \textquotedblleft symmetrized action form\textquotedblright : 
\begin{equation*}
\beta _{H}=\frac{1}{2}(pdx-xdp)-Hdt
\end{equation*}%
We claim (somewhat speculatively...) that had Schr\"{o}dinger used $\beta
_{H}$ instead of $\alpha _{H}$ he could very well have landed, not with the
equation (\ref{erwin1}), but rather with the phase-space equation (\ref%
{erwin2}), and could hence have led him to deformation quantization!

Let us justify our claims from a rigorous mathematical point of view.

\subsection*{Notations}

We denote by $\sigma $ the canonical symplectic form on the phase space $%
\mathbb{R}_{z}^{2n}=\mathbb{R}_{x}^{n}\times \mathbb{R}_{p}^{n}$:%
\begin{equation*}
\sigma (z,z^{\prime })=px^{\prime }-p^{\prime }x\text{ \ if \ }z=(x,p)\text{%
, }z^{\prime }=(x^{\prime }p^{\prime })
\end{equation*}%
where $x=(x_{1},...,x_{1})$, $p=(p_{1},...,p_{n})$; we are using the
\textquotedblleft dotless dot-product\textquotedblright\ notation $%
xp=x_{1}p_{1}+\cdot \cdot \cdot +x_{n}p_{n}$. The generalized gradients $%
\partial _{x}$ and $\partial _{p}$ are defined by $\partial _{x}=(\partial
/\partial x_{1},...,\partial /\partial x_{n})$ and $\partial _{p}=(\partial
/\partial p_{1},...,\partial /\partial p_{n})$.

We denote by $Sp(n)$ the real symplectic group; it consists of all linear
automorphisms $S$ of $\mathbb{R}_{z}^{2n}$ such that $\sigma (Sz,Sz^{\prime
})=\sigma (z,z^{\prime })$ for all $z,z^{\prime }$.

$\mathcal{S}(\mathbb{R}^{m})$ is the Schwartz space of rapidly decreasing
functions on $\mathbb{R}^{m}$; its dual $\mathcal{S}^{\prime }(\mathbb{R}%
^{m})$ is the space of tempered distributions. Functions on $\mathbb{R}%
_{x}^{n}$ or $\mathbb{R}_{p}^{n}$ will be denoted by small Greek letters $%
\psi ,\phi ,...$ while functions on $\mathbb{R}_{z}^{2n}$ will be denoted by
capital Greek letters, e.g., $\Psi $.

For the notions of Weyl calculus that are being used here, see Folland \cite%
{Folland} or \cite{Wong}; we are using the notations and normalizations of
Littlejohn \cite{Littlejohn}. For a review of deformation quantization see
the preface in Zachos \textit{et al}. \cite{ZFC}. 

\section{Phase-Space Representation of \textbf{H}$_{n}$}

Recall that one of the modern ways to justify the Schr\"{o}dinger
quantization rules $x_{j}\longrightarrow x_{j},$ $p_{j}\longrightarrow
-i\hbar (\partial /\partial x_{j})$ is to construct the Schr\"{o}dinger
representation of the Heisenberg group $\mathbf{H}_{n},$ that is $\mathbb{R}%
_{z}^{2n}\times \mathbb{R}_{t}$ equipped with the group law 
\begin{equation}
(z,t)\cdot (z^{\prime },t)=(z+z^{\prime },t+t^{\prime }+\tfrac{1}{2}\sigma
(z,z^{\prime }))\text{.}  \label{hglaw1}
\end{equation}%
One proceeds as follows: consider the \textquotedblleft translation
Hamiltonian\textquotedblright\ $H_{z_{0}}=\sigma (z,z_{0})$; the flow it
determines are the translations $T(tz_{0}):z\longmapsto z+tz_{0}$; they act
on functions defined on $\mathbb{R}_{z}^{2n}$ by the rule%
\begin{equation*}
T(tz_{0})\Psi _{0}(z)=\Psi _{0}(z-tz_{0})\text{.}
\end{equation*}%
In (traditional) quantum mechanics Hilbert spaces and phases play a crucial
role; one \textquotedblleft quantizes\textquotedblright\ the operators $%
T(tz_{0})$ by letting them act on $\psi _{0}\in L^{2}(\mathbb{R}_{x}^{n})$
via the Heisenberg--Weyl operators $\widehat{T}(z_{0})$ defined by%
\begin{equation*}
\widehat{T}(tz_{0})\psi _{0}(x)=e^{\frac{i}{\hbar }\varphi
(z,t)}T(tz_{0})\psi _{0}(x);
\end{equation*}%
here $\varphi (z,t)$ is the increase in action when one goes straight from
the point $z-tz_{0}$ to the point $z$, that is%
\begin{equation}
\varphi (z,t)=\int_{-t}^{0}pdx-H_{z_{0}}dt=p_{0}xt-\frac{t^{2}}{2}p_{0}x_{0}%
\text{;}  \label{fihab}
\end{equation}%
thus%
\begin{equation}
\widehat{T}(tz_{0})\psi _{0}(x)=e^{\frac{i}{\hbar }(p_{0}xt-\frac{t^{2}}{2}%
p_{0}x_{0})}\psi _{0}(x-tx_{0})\text{.}  \label{terwin}
\end{equation}%
The Schr\"{o}dinger representation of $\mathbf{H}_{n}$ in $L^{2}(\mathbb{R}%
_{x}^{n})$ is the mapping%
\begin{equation*}
T_{\text{Sch}}:\mathbf{H}_{n}\longrightarrow \mathcal{U}(L^{2}(\mathbb{R}%
_{x}^{n}\mathcal{))}
\end{equation*}%
($\mathcal{U}(L^{2}(\mathbb{R}_{x}^{n}\mathcal{))}$ the unitary operators on 
$L^{2}(\mathbb{R}_{x}^{n}\mathcal{)}$) defined by%
\begin{equation}
T_{\text{Sch}}(z_{0},t_{0})\psi _{0}(x)=e^{\frac{i}{\hbar }t_{0}}\widehat{T}%
(z_{0})\psi _{0}(x);  \label{repsch}
\end{equation}%
one proves that $T_{\text{Sch}}$ is a unitary and irreducible
representation; a famous theorem of Stone and von Neumann (see \cite%
{Folland,Wallach} for a proof) asserts that it is, up to unitary
equivalences, the \textit{only} irreducible representation of $\mathbf{H}%
_{n} $ in $L^{2}(\mathbb{R}_{x}^{n}\mathcal{)}$. But this theorem does not
prevent us from constructing non-trivial irreducible representations of $%
\mathbf{H}_{n}$ in \textit{other} Hilbert spaces; we will come back to this
essential point in a moment, but let us first note that Schr\"{o}dinger's
equation for the displacement Hamiltonian $H_{z_{0}}=\sigma (z,z_{0})$, and
hence the quantum rules%
\begin{equation}
x_{j}\longrightarrow x_{j}\text{ \ \ \ , \ \ \ }p_{j}\longrightarrow -i\hbar 
\frac{\partial }{\partial x_{j}}  \label{erwinusual}
\end{equation}%
now follow from formula (\ref{terwin}): an immediate calculation shows that
the function $\psi (x,t)=\widehat{T}(tz_{0})\psi _{0}(x)$ is a solution of 
\begin{equation*}
i\hbar \frac{\partial \psi }{\partial t}=H_{z_{0}}\left( x,-i\hbar \partial
_{x}\right) \psi \text{ \ , \ }\psi (x,0)=\psi _{0}(x).
\end{equation*}

Let us quantize the translation operators $T(tz_{0})$ in a different way. We 
\emph{redefine} $\widehat{T}(tz_{0})$ by letting it act, not on $L^{2}(%
\mathbb{R}_{x}^{n}\mathcal{)}$, but on $L^{2}(\mathbb{R}_{z}^{2n}\mathcal{)}$%
, by the formula 
\begin{equation*}
\widehat{T}_{\text{ph}}(tz_{0})\Psi _{0}(z)=e^{\frac{i}{\hbar }\varphi
^{\prime }(z,t)}T(tz_{0})\Psi _{0}(z)
\end{equation*}%
(the subscript \textquotedblleft ph\textquotedblright\ stands for
\textquotedblleft \textit{ph}ase space\textquotedblright ), and replacing
the phase (\ref{fihab}) by integrating, not the Poincar\'{e}--Cartan form $%
\alpha _{H_{z_{0}}}$ but its symmetrized variant 
\begin{equation*}
\beta _{H_{z_{0}}}=\frac{1}{2}(pdx-xdp)-H_{z_{0}}dt\text{.}
\end{equation*}%
This yields after a trivial calculation%
\begin{equation}
\varphi ^{\prime }(z,t)=-\frac{1}{2}H_{z_{0}}(z)t=-\frac{1}{2}\sigma
(z,z_{0})t.  \label{finew}
\end{equation}%
Summarizing, we have defined 
\begin{equation}
\widehat{T}_{\text{ph}}(tz_{0})\Psi _{0}(z)=e^{-\frac{i}{2\hbar }\sigma
(z,z_{0})t}\Psi _{0}(z-tz_{0}).  \label{hwnew}
\end{equation}%
What partial differential equation does the function $\Psi =\widehat{T}_{%
\text{ph}}(tz_{0})\Psi _{0}$ satisfy? Performing a few calculations one
checks that it satisfies the multi-dimensional analogue of the phase-space
Schr\"{o}dinger equation (\ref{erwin2}) of the introduction, namely%
\begin{equation}
i\hbar \frac{\partial \Psi }{\partial t}=H\left( \tfrac{x}{2}+i\hbar
\partial _{p},\tfrac{p}{2}-i\hbar \partial _{x}\right) \Psi .  \label{erwin3}
\end{equation}

We are going to prove the following:

\begin{description}
\item[(A)] The operators $\widehat{T}_{\text{ph}}(tz_{0})$ correspond to a
new irreducible unitary representation of the Heisenberg group $\mathbf{H}%
_{n}$ on a closed subspace of $L^{2}(\mathbb{R}_{z}^{2n}\mathcal{)}$ (which
is unitarily equivalent to the Schr\"{o}dinger representation via Stone--von
Neumann's theorem).

\item[(B)] The phase-space Schr\"{o}dinger equation (\ref{erwin3}) is
closely related to deformation quantization, in fact to an extension of the
usual Weyl calculus on $L^{2}(\mathbb{R}_{x}^{n}\mathcal{)}$ to $L^{2}(%
\mathbb{R}_{z}^{2n}\mathcal{)}$, for which the operators $H\left( \frac{x}{2}%
+i\hbar \partial _{p},\frac{p}{2}-i\hbar \partial _{x}\right) $ are
perfectly well-defined.
\end{description}

\section{The Irreducible Unitary Representation $\widehat{T}_{\text{ph}}$}

We define the phase-space representation of $\mathbf{H}_{n}$ in analogy with
(\ref{repsch}) by 
\begin{equation}
\widehat{T}_{\text{ph}}(z_{0},t_{0})\Psi _{0}(z)=e^{\frac{i}{\hbar }t_{0}}%
\widehat{T}_{\text{ph}}(tz_{0})\Psi _{0}(z)\text{.}  \label{tps}
\end{equation}%
Clearly $\widehat{T}_{\text{ph}}(z_{0},t_{0})$ is a unitary operator;
moreover a straightforward calculation shows that%
\begin{equation*}
\widehat{T}_{\text{ph}}(z_{0},t_{0})\widehat{T}_{\text{ph}}(z_{1},t_{1})=e^{%
\frac{i}{2\hbar }\sigma (z_{0},z_{1})}\widehat{T}_{\text{ph}%
}(z_{0}+z_{1},t_{0}+t_{1}+\tfrac{1}{2}\sigma (z_{0},z_{1}))
\end{equation*}%
so that $\widehat{T}_{\text{ph}}$ is indeed a representation of $\mathbf{H}%
_{n}$ on some subspace of $L^{2}(\mathbb{R}_{z}^{2n})$. We are going to show
that this representation is unitarily equivalent to the Schr\"{o}dinger
representation, and hence irreducible.

Let $\phi \in \mathcal{S}(\mathbb{R}_{x}^{n})$ be normalized: $||\phi
||_{L^{2}(\mathbb{R}_{x}^{n})}^{2}=1$. To $\phi $ we associate the operator $%
V_{\phi }:L^{2}(\mathbb{R}_{x}^{n}\mathcal{)}\longrightarrow L^{2}(\mathbb{R}%
_{z}^{2n}\mathcal{)}$ defined by%
\begin{equation*}
V_{\phi }\psi (z)=\left( \tfrac{\pi \hbar }{2}\right) ^{n/2}W(\psi ,%
\overline{\phi })(\tfrac{1}{2}z)
\end{equation*}%
where $W(\psi ,\overline{\phi })$ is the Wigner--Moyal function (Folland 
\cite{Folland}):%
\begin{equation*}
W(\psi ,\overline{\phi })(x,p)=\left( \tfrac{1}{2\pi \hbar }\right) ^{n}\int
e^{-\tfrac{i}{\hbar }\left\langle p,y\right\rangle }\psi (x+\tfrac{1}{2}%
y)\phi (x-\tfrac{1}{2}y)d^{n}y\text{.}
\end{equation*}%
It turns out that $V_{\phi }$ is an extension of the \textquotedblleft
coherent-state representation\textquotedblright\ to which it reduces,up to
the factor $\exp (-ipx/\hbar )$ if one takes for $\phi $ a the Gaussian%
\begin{equation}
\phi _{0}(x)=\left( \tfrac{1}{\pi \hbar }\right) ^{n/4}e^{-\frac{1}{2\hbar }%
|x|^{2}}\text{.}  \label{fizero}
\end{equation}%
In fact, a straightforward calculation shows that%
\begin{equation}
V_{\phi }\psi (z)=e^{-\frac{i}{2\hbar }px}U_{\phi }\psi (z)  \label{vif}
\end{equation}%
where $U_{\phi }$ is 
\begin{equation}
U_{\phi }\psi (z)=\left( \tfrac{1}{2\pi \hbar }\right) ^{n/2}\int e^{\frac{i%
}{\hbar }\left\langle p,x-x^{\prime }\right\rangle }\phi (x-x^{\prime })\psi
(x^{\prime })d^{n}x^{\prime }\text{.}  \label{ufi}
\end{equation}%
It follows from the properties of $U_{\phi }$ (see for instance \cite{Naza})
that

\begin{enumerate}
\item The transform $V_{\phi }$ is an isometry:%
\begin{equation}
(V_{\phi }\psi ,V_{\phi }\psi ^{\prime })_{L^{2}(\mathbb{R}_{z}^{2n})}=(\psi
,\psi ^{\prime })_{L^{2}(\mathbb{R}_{x}^{n})}  \label{Parseval}
\end{equation}%
holds for all $\psi ,\psi ^{\prime }\in \mathcal{S}(\mathbb{R}_{x}^{n})$

\item $V_{\phi }$ extends into an isometric operator $L^{2}(\mathbb{R}%
_{x}^{n})\longrightarrow L^{2}(\mathbb{R}_{z}^{2n})$ and 
\begin{equation}
V_{\phi }^{\ast }V_{\phi }=I\text{ \ on \ }L^{2}(\mathbb{R}_{x}^{n});
\label{parsevalin}
\end{equation}

\item The range $\mathcal{H}_{\phi }$ of $V_{\phi }$ is closed in $L^{2}(%
\mathbb{R}_{z}^{2n})$ (and is hence a Hilbert space), and $P=V_{\phi
}V_{\phi }^{\ast }$ is the orthogonal projection on the Hilbert space $%
\mathcal{H}_{\phi }.$
\end{enumerate}

To show that $\widehat{T}_{\text{ph}}$ is unitarily equivalent to the Schr%
\"{o}dinger representation, it suffices to show that the operators $\widehat{%
T}_{\text{ph}}(z_{0})=\widehat{T}_{\text{ph}}(z_{0},0)$ and $T_{\text{Sch}%
}(z_{0})=T_{\text{Sch}}(z_{0},0)$ are such that 
\begin{equation}
\widehat{T}_{\text{ph}}(z_{0})V_{\phi }=V_{\phi }\widehat{T}_{\text{Sch}%
}(z_{0})  \label{unit}
\end{equation}%
for every $z_{0}$. Now,%
\begin{eqnarray*}
\widehat{T}_{\text{ph}}(z_{0})V_{\phi }\psi (z) &=&e^{-\frac{i}{2\hbar }%
\sigma (z,z_{0})}e^{-\frac{i}{2\hbar }px}U_{\phi }\psi (z-z_{0}) \\
&=&\left( \tfrac{1}{2\pi \hbar }\right) ^{n/2}e^{-\frac{i}{\hbar }(p_{0}x-%
\frac{1}{2}p_{0}x_{0})}\times \\
&&\int e^{\frac{i}{\hbar }(p-p_{0})(x-x_{0}-x^{\prime })}\phi
(x-x_{0}-x^{\prime })\psi (x^{\prime })d^{n}x^{\prime }
\end{eqnarray*}%
and setting $x^{\prime \prime }=x^{\prime }+x_{0}$ in the integral this is%
\begin{eqnarray*}
\widehat{T}_{\text{ph}}(z_{0})V_{\phi }\psi (z) &=&\left( \tfrac{1}{2\pi
\hbar }\right) ^{n/2}e^{-\frac{i}{\hbar }(p_{0}x-\frac{1}{2}%
p_{0}x_{0})}\times \\
&&\int e^{\frac{i}{\hbar }(p-p_{0})(x-x^{\prime \prime })}\phi (x-x^{\prime
\prime })\psi (x^{\prime \prime }-x_{0})d^{n}x^{\prime }
\end{eqnarray*}%
hence%
\begin{equation*}
\widehat{T}_{\text{ph}}(z_{0})(V_{\phi }\psi )(z)=V_{\phi }(T_{\text{Sch}%
}(z_{0})\psi )(z)
\end{equation*}%
which was to be proven.

\begin{remark}
The Hilbert space $\mathcal{H}_{\phi }$ is smaller than $L^{2}(\mathbb{R}%
_{z}^{2n})$; for instance if we chose for $\phi $ the Gaussian (\ref{fizero}%
) then one proves \cite{Naza} that the range of the transform $U_{\phi }$
defined by (\ref{ufi}) consists of all $\Psi \in L^{2}(\mathbb{R}_{z}^{2n})$
such that $\exp (p^{2}/2\hbar )$ is anti-analytic. It follows that $\mathcal{%
H}_{\phi _{0}}$ which is the range of $V_{\phi }=\exp (-ipx/2\hbar )U_{\phi }
$ consists of all $\Psi \in L^{2}(\mathbb{R}_{z}^{2n})$ such that%
\begin{equation*}
\frac{\partial }{\partial z_{j}}(e^{\frac{1}{2\hbar }|z|^{2}}\Psi (z))=0%
\text{ \ , \ }1\leq j\leq n\text{.}
\end{equation*}
\end{remark}

Moreover, a few calculations, using for instance (\ref{vif}) and (\ref{ufi})
show that we have%
\begin{equation}
\left( \tfrac{x}{2}+i\hbar \partial _{p}\right) V_{\phi }\psi =V_{\phi
}(x\psi )\text{ \ , \ }\left( \tfrac{p}{2}-i\hbar \partial _{x}\right)
V_{\phi }\psi =V_{\phi }(-i\hbar \partial _{x}\psi );  \label{erwin4}
\end{equation}%
the transform $V_{\phi }$ thus takes the usual quantization rules (\ref%
{erwinusual}) to the phase-space quantization rules.%
\begin{equation*}
x\longrightarrow \frac{x}{2}+i\hbar \partial _{p}\text{ \ , \ }%
x\longrightarrow \frac{p}{2}-i\hbar \partial _{x}\text{.}
\end{equation*}

\section{Extended Weyl Calculus}

In standard Weyl calculus one associates to a \textquotedblleft
symbol\textquotedblright\ $a$ having some some suitable growth properties
for $p\rightarrow \infty $ class a pseudo-differential operator 
\begin{equation*}
\widehat{A}=a^{w}:\mathcal{S}(\mathbb{R}_{x}^{n})\longrightarrow \mathcal{S}(%
\mathbb{R}_{x}^{n})
\end{equation*}%
defined by the kernel 
\begin{equation*}
K_{\widehat{A}}(x,y)=\left( \tfrac{1}{2\pi \hbar }\right) ^{N/2}\int e^{%
\frac{i}{\hbar }p(x,y)}a(\tfrac{1}{2}(x+y),p)d^{N}p.
\end{equation*}%
One proves that 
\begin{equation}
\widehat{A}\psi (x)=\left( \tfrac{1}{2\pi \hbar }\right) ^{n}\dint \tilde{a}%
(z_{0})\widehat{T}_{\text{Sch}}(z_{0})\psi (x)d^{2n}z_{0}  \label{weyl1}
\end{equation}%
for $\psi \in \mathcal{S}(\mathbb{R}_{x}^{n})$ (the integral being
interpreted as an \textquotedblleft oscillatory integral\textquotedblright ,
see \textit{e.g.} \cite{Folland,Wong}). In formula (\ref{weyl1}) $\tilde{a}$
(the \textquotedblleft twisted\textquotedblright\ Weyl symbol) is the
symplectic Fourier-transform of $a$:%
\begin{equation}
\tilde{a}(z)=\mathcal{F}_{\sigma }a(z)=\left( \tfrac{1}{2\pi \hbar }\right)
^{n}\dint e^{-\frac{i}{\hbar }\sigma (z,z^{\prime })}a(z^{\prime
})d^{2n}z^{\prime }  \label{sft}
\end{equation}%
and $\widehat{T}_{\text{Sch}}(z_{0})=\widehat{T}_{\text{Sch}}(z_{0},0)$ is
the Heisenberg--Weyl operator (\ref{terwin}).

The discussion above suggests that we might now be able to make $\widehat{A}$
to act, not only on functions of $x,$ but also on functions $\Psi \in 
\mathcal{S(}\mathbb{R}_{z}^{2n})$ by defining 
\begin{equation}
\widehat{A}_{\text{Sch}}\Psi (z)=\left( \tfrac{1}{2\pi \hbar }\right)
^{n}\dint \tilde{a}(z_{0})\widehat{T}_{\text{Sch}}(z_{0})\Psi (z)d^{2n}z_{0}
\label{weyl2}
\end{equation}%
where we have set%
\begin{equation*}
\widehat{T}_{\text{Sch}}(z_{0})\Psi (z)=e^{\frac{i}{\hbar }(p_{0}x-\frac{1}{2%
}p_{0}x_{0})}\Psi (z-z_{0})\text{.}
\end{equation*}%
It turns out that it is better for our purposes to define instead the
operator%
\begin{equation}
\widehat{A}_{\text{ph}}\Psi (z)=\left( \tfrac{1}{2\pi \hbar }\right)
^{n}\dint \tilde{a}(z_{0})\widehat{T}_{\text{ph}}(z_{0})\Psi (z)d^{2n}z_{0}%
\text{.}  \label{weyl3}
\end{equation}%
that is%
\begin{equation}
\widehat{A}_{\text{ph}}\Psi (z)=\left( \tfrac{1}{2\pi \hbar }\right)
^{n}\dint e^{-\frac{i}{2\hbar }\sigma (z,z_{0})}\mathcal{F}_{\sigma
}a(z_{0})\Psi (z-z_{0})d^{2n}z_{0}\text{.}  \label{weyl4}
\end{equation}

It turns out that this formula ids the fundamental link between the theory
sketched above with deformatio quantization.

\section{Relation With Deformation Quantization}

Recall now that if $\widehat{A}=a^{w}$ and $\widehat{B}=b^{w}$ are the Weyl
operators with symbols $a$ and $b$, respectively, then the twisted symbol $%
\tilde{c}=\mathcal{F}_{\sigma }c$ of the compose $\widehat{C}=\widehat{A}%
\circ \widehat{B}$ is given by 
\begin{equation*}
\tilde{c}(z)=\left( \tfrac{1}{2\pi \hbar }\right) ^{n}\dint e^{-\frac{i}{%
2\hbar }\sigma (z,z^{\prime })}a(z-z^{\prime })b(z^{\prime })d^{2n}z^{\prime
}\text{;}
\end{equation*}%
since $\mathcal{F}_{\sigma }$ is an involution, we have $c=\mathcal{F}%
_{\sigma }\tilde{c}$ and one verifies that 
\begin{equation*}
c(z)=\left( \tfrac{1}{4\pi \hbar }\right) ^{2n}\dint e^{\frac{i}{2\hbar }%
\sigma (z^{\prime },z^{\prime \prime })}a(z+\tfrac{1}{2}z^{\prime })b(z-%
\tfrac{1}{2}z^{\prime \prime })d^{2n}z^{\prime }d^{2n}z^{\prime \prime }%
\text{;}
\end{equation*}%
using expansions in Taylor series and repeated integrations by parts this
can be rewritten in terms of the \textquotedblleft Janus
operator\textquotedblright\ $\overleftarrow{\partial _{x}}\overrightarrow{%
\partial _{p}}-\overleftarrow{\partial _{p}}\overrightarrow{\partial _{x}}$
as 
\begin{equation*}
c(z)=a(z)\exp \left[ \tfrac{i\hbar }{2}(\overleftarrow{\partial _{x}}\cdot 
\overrightarrow{\partial _{p}}-\overleftarrow{\partial _{p}}\cdot 
\overrightarrow{\partial _{x}}\right] b(z)\text{;}
\end{equation*}%
in deformation quantization this is called the star-product (or Moyal
product) $a\bigstar b$ of the symbols $a$ and $b$. Thus formula (\ref{weyl4}%
) says that our extended Weyl calculus can be expressed in terms of the
star-product in the following very simple way:%
\begin{equation}
\widehat{A}_{\text{ph}}\Psi =\mathcal{F}_{\sigma }(a\bigstar \Psi ).
\label{Weyl5}
\end{equation}

\section{Discussion and Concluding Remarks}

The relevance of Torres-Vega and Frederick's Schr\"{o}dinger equation in
phase space is justified not only because it is consistent with Stone and
von Neumann's theorem on the irreducible representations of the Heisenberg
group, but also because it is a variant of deformation quantization; any
advance in one of these theory will thus lead to an advance in the other. 

\begin{acknowledgement}
This work has been partially supported by a grant of the Max-Planck-Institut
f\"{u}r Gravitationsphysik (Albert-Einstein-Institut, Golm). I wish to thank
Prof. Hermann Nicolai for his kind hospitality.
\end{acknowledgement}

\end{document}